\documentclass[11pt,a4paper,reqno]{amsart}

\usepackage{amsfonts}
\usepackage{mathrsfs}
\usepackage{graphicx}
\usepackage{amsmath} 
\usepackage{graphics}
\usepackage{epsfig}
\usepackage{anysize}
\usepackage{lipsum}
\usepackage{wrapfig}
\usepackage{color}

\marginsize{3cm}{3cm}{3cm}{3cm}
\newtheorem{theorem}{Theorem}[section]
\newtheorem{proposition}{Proposition}[section]

\newtheorem{definition}{Definition}[section]

\makeatletter

\begin{document}

\setcounter{page}{1}
\vspace{2cm}
\author[\hspace{0.7cm}\centerline{TWMS J. App. Eng.
Math. V.xx, N.xx, 20xx}]{V. Chinnadurai$^1$, A. Arulselvam$^2$,  \S }
\title[\centerline{V. Chinnadurai, A. Arulselvam: Rough cubic Pythagorean fuzzy sets in semigroup.
\hspace{0.5cm}}]{Rough cubic Pythagorean fuzzy sets in semigroup}

\thanks{\noindent $^1$ V. Chinnadurai, Professor, Department of Mathematics, Annamalai University, India.\\
\indent \,\,\, e-mail: chinnaduraiau@gmail.com; ORCID: https://orcid.org/0000-0002-6047-6348. \\
\indent $^2$ A. Arulselvam, Ph.D Research scholar, Department of Mathematics, Annamalai University, India.\\
\indent \,\,\, e-mail: arulselvam.a91@gmail.com; ORCID: https://orcid.org/0000-0001-8383-2889. of the second \\
\indent \,\,\, author\\
\indent \S \, Manuscript received: Month Day, Year; accepted: Month Day, Year. \\
\indent \,\,\, TWMS Journal of Applied and Engineering Mathematics, Vol.xx,
No.xx; \copyright\ I\c s\i k University, Depart-\\ 
\indent \,\,\,ment of Mathematics, 20xx; all rights reserved.}

%
\begin{abstract}
In this paper, we intend the concept of rough cubic Pythagorean fuzzy ideals in the semigroup. By using this notion, we discuss lower approximation and upper approximation of cubic Pythagorean fuzzy left (right) ideals, bi-ideals, interior ideals, and study some of their related properties in detail.
 
\end{abstract}
\maketitle 
\bigskip

%
\section{Introduction}
In 1965, Zadeh\cite{Z1,Z2} introduced the concept of fuzzy sets, and later in 1975 developed the interval-valued fuzzy set an extension of a fuzzy set. A semigroup is an algebraic structure consisting of non-empty sets together with an associative binary operation. The formal study of semigroup began in the early twentieth century. Pawlak\cite{P} initiated the fundamental concept of rough set in 1982. Dubois and Prade\cite{DP} developed the concepts of rough fuzzy sets based on Pawlak approximation in 1989. In 1997 the idea of rough ideals in semigroup was presented by Kuroki\cite{K}. Jun et.al\cite{JKY} initiated the new notion called a cubic set, which is a combination of interval-valued fuzzy set and fuzzy set and discussed some of its related properties in 2012. Jun and Khan\cite{JK} introduced the notion of cubic ideals in semigroup in 2013. Yager\cite{Y} initiated the notion of Pythagorean fuzzy set in 2013 the concept of sum of the squares of membership and non-membership belongs to the unit interval [0,1]. Hussain et.al\cite{HMI} initiated the notions of rough Pythagorean fuzzy ideals in semigroup in 2019. Garg\cite{G1,G2,G3} exploited Pythagorean fuzzy sets to solve multi-criteria decision making problems. Zhang and Xu\cite{ZX} presented TOPSIS method in Pythagorean fuzzy sets to rank the alternatives.\\
The main aim of this paper is to study the notion of rough cubic Pythagorean fuzzy sets in the semigroup and investigate some of its related properties. We study the properties of rough cubic Pythagorean fuzzy sets in semigroup by using congruence relation. Also, we prove some interesting properties of the rough cubic Pythagorean fuzzy left (right) ideal, ideal, bi-ideal, and interior-ideal. 
\section{Preliminaries} 
The basic concepts of Rough set(RS), Pythagorean fuzzy set(PFS), Cubic Pythagorean fuzzy set(CPFS), Rough Pythagorean fuzzy set(RCPFS) are referred \cite{P}, \cite{Y}, \cite{AAASH}, \cite{HMI}  are respectively.\\
\begin{definition}\cite{A}
Let $X$ be a universe of discourse, An \textbf{intuitionistic fuzzy set}(IFS) $A$ in $X$ is an object having the form. $A=\left\{z,\zeta_{A}(z),\eta_{A}(z)/ z\in X\right\}$. where the mapping $\zeta:X\rightarrow [0,1]$ and $\eta:X\rightarrow [0,1]$ represent the degree of membership and non-membership of the object $z\in X$ to the set A respectively with the condition $0\leq\zeta_{A}(z)+\eta_{A}(z)\leq 1$. for all $z\in X$ for the sake of simplicity an IFS is denoted by $A=(\zeta_{A}(z),\eta_{A}(z))$.
\end{definition}
\begin{definition}
Let $\widetilde{P}=\left(\zeta_{\widetilde{p}},\eta_{\widetilde{p}}\right)=\left\{\left\langle z,\zeta_{\widetilde{p}}(z),\eta_{\widetilde{p}}(z)\right\rangle/z\in S\right\}$ be the Interval Valued Pythagorean fuzzy set(IVPFS) of S, where $\zeta_{\widetilde{p}}(z)=\zeta_{\widetilde{p}^{-}}(z),\zeta_{\widetilde{p}^{+}}(z)$ and $\eta_{\widetilde{p}}(z)=\eta_{\widetilde{p}^{-}}(z),\eta_{\widetilde{p}^{+}}(z)$. Then RIVPFS of S is denoted as $App\left(\widetilde{p}\right)=\left(\underline{App}(\widetilde{p}),\overline{App}(\widetilde{p})\right)$ the lower approximation is defined as $\underline{App}(\widetilde{p})=\left\{\left\langle z,\underline{\zeta_{\widetilde{p}}(z)},\underline{\eta_{\widetilde{p}}(z)}\right\rangle/z\in S\right\}$ where $\underline{\zeta}_{\widetilde{p}}=\bigwedge\limits_{z^{'}\in [z]_{\omega}}\zeta_{\widetilde{p}}(z^{'})$ and $\underline{\eta}_{\widetilde{p}}=\bigvee\limits_{z^{'}\in [z]_{\omega}}\eta_{\widetilde{p}}(z^{'})$ with condition that $0\leq\underline{\left(\zeta_{\widetilde{p}}(z)\right)}^{2}+\underline{\left(\eta_{\widetilde{p}}(z)\right)}^{2}\leq 1$ and the lower approximation is defined as $\overline{App}(\widetilde{p})=\left\{\left\langle z,\overline{\zeta_{\widetilde{p}}(z)},\overline{\eta_{\widetilde{p}}(z)}\right\rangle/z\in S\right\}$ where $\overline{\zeta}_{\widetilde{p}}=\bigvee\limits_{z^{'}\in [z]_{\omega}}\zeta_{\widetilde{p}}(z^{'})$ and $\overline{\eta}_{\widetilde{p}}=\bigwedge\limits_{z^{'}\in [z]_{\omega}}\eta_{\widetilde{p}}(z^{'})$ with condition that $0\leq\overline{\left(\zeta_{\widetilde{p}}(z)\right)}^{2}+\overline{\left(\eta_{\widetilde{p}}(z)\right)}^{2}\leq 1$.
\end{definition}
Throughout this paper $S$ denotes the semigroup.
\begin{definition}
Let $P^{\Box}_{1}=\left(\zeta_{p^{\Box}_{1}},\eta_{p^{\Box}_{1}}\right)$ and $P^{\Box}_{2}=\left(\zeta_{p^{\Box}_{2}},\eta_{p^{\Box}_{2}}\right)$ be any two CPFS on S. Then, the composition of $P^{\Box}_{1}$ and $P^{\Box}_{2}$ is defined as $P^{\Box}_{1}\circ P^{\Box}_{2}=\left(\zeta_{p^{\Box}_{1}}\circ\zeta_{p^{\Box}_{2}},\eta_{p^{\Box}_{1}}\circ\eta_{p^{\Box}_{2}}\right)$ where, $\left(\zeta_{p^{\Box}_{1}}\circ\zeta_{p^{\Box}_{2}}\right)(z)=\bigvee\limits_{z=z_{1}z_{2}}\left[\zeta_{p^{\Box}_{1}}(z_{1})\bigwedge\zeta_{p^{\Box}_{2}}(z_{2})\right]$\\
$\left(\eta_{p^{\Box}_{1}}\circ\eta_{p^{\Box}_{2}}\right)(z)=\bigwedge\limits_{z=z_{1}z_{2}}\left[\eta_{p^{\Box}_{1}}(z_{1})\bigvee\eta_{p^{\Box}_{2}}(z_{2})\right]$.
\end{definition}
\section{Rough cubic pythagorean fuzzy sets (RCPFS) in semigroup}
A equivalence relation $\omega$ on $S$ is said to be a congruence relation denoted as ${CR}_{\omega}$ if for all $x,z_{1},z_{2}\in S$ such that $z_{1},z_{2}\in \omega$ $\Rightarrow z_{1}x,z_{2}x\in\omega$ and $xz_{1},xz_{2}\in\omega$ The congruence class of an object $z\in S$ is denoted by $[z]_{\omega}$. For a $CR_{\omega}$ on $S$, we have $[z_{1}]_{\omega}[z_{2}]_{\omega}\subseteq[z_{1}z_{2}]_{\omega}$ and the $CR_{\omega}$ on $S$ is called complete if $[z_{1}]_{\omega}[z_{2}]_{\omega}=[z_{1}z_{2}]_{\omega}$. $\forall z_{1},z_{2}\in S$.
\begin{definition}
Let $P^{\Box}=\left(\zeta_{p^{\Box}},\eta_{p^{\Box}}\right)=\left\{\left\langle 
z_{1},\left[\zeta_{\widetilde{p}}(z_{1}),\eta_{\widetilde{p}}(z_{1})\right],\left(\zeta_{p}(z_{1}),\eta_{p}(z_{1})\right)\right\rangle/z_{1}\in S\right\}$ be the CPFS in $S$, where $\zeta_{\widetilde{p}}(z_{1})=\left(\zeta^{-}_{p}(z_{1}),\zeta^{+}_{p}(z_{1})\right)$ and $\eta_{\widetilde{p}}(z_{1})=\left(\eta^{-}_{p}(z_{1}),\eta^{+}_{p}(z_{1})\right)$. Then a RCPFS on $S$ is denoted by $App(P^{\Box})=\left(\underline{App}(P^{\Box}),\overline{App}(P^{\Box})\right)$. The lower approximation is defined as $\underline{App}(P^{\Box})=\left(\underline{\zeta_{p^{\Box}}},\underline{\eta_{p^{\Box}}}\right)=\left\{\left\langle 
z_{1},\left[\underline{\zeta_{\widetilde{p}}}(z_{1}),\underline{\eta_{\widetilde{p}}}(z_{1})\right],\left(\underline{\zeta_{p}}(z_{1}),\underline{\eta_{p}}(z_{1})\right)\right\rangle/z_{1}\in S\right\}$ where $\underline{\zeta_{\widetilde{p}}}(z)=\bigwedge\limits_{z^{'}\in[z]_{\omega}}\zeta_{\widetilde{p}}(z^{'})$ and $\underline{\eta_{\widetilde{p}}}(z)=\bigvee\limits_{z^{'}\in[z]_{\omega}}\eta_{\widetilde{p}}(z^{'})$
$\underline{\zeta_{p}}(z)=\bigwedge\limits_{z^{'}\in[z]_{\omega}}\zeta_{p}(z^{'})$ and $\underline{\eta_{p}}(z)=\bigvee\limits_{z^{'}\in[z]_{\omega}}\eta_{p}(z^{'})$ with the condition $0\leq\left(\underline{\zeta_{p}}(z)\right)^{2}+\left(\underline{\eta_{p}}(z)\right)^{2}\leq 1$ and the upper approximation is defined as $\overline{App}(P)=\left\{\left\langle z,\overline{\zeta_{p}}(z),\overline{\eta_{p}}(z)\right\rangle/z\in S\right\}$. where $\overline{\zeta_{\widetilde{p}}}(z)=\bigvee\limits_{z^{'}\in[z]_{\omega}}\zeta_{\widetilde{p}}(z^{'})$ and $\overline{\eta_{\widetilde{p}}}(z)=\bigwedge\limits_{z^{'}\in[z]_{\omega}}\eta_{\widetilde{p}}(z^{'})$
 $\overline{\zeta_{p}}(z)=\bigvee\limits_{z^{'}\in[z]_{\omega}}\zeta_{p}(z^{'})$ and $\overline{\eta_{p}}(z)=\bigwedge\limits_{z^{'}\in[z]_{\omega}}\eta_{p}(z^{'})$ with the condition $0\leq\left(\overline{\zeta_{p}}(z)\right)^{2}+\left(\overline{\eta_{p}}(z)\right)^{2}\leq 1$.
\end{definition}
\begin{proposition}
The lower approximation and upper approximation of the CPFS $P^{\Box}$ on $S$ are CPFS of a quotient set $S/\omega$
\begin{proof}
The membership and non-membership grades of lower approximation i.e., $\underline{App}(P^{\Box})$ from definition 3.1 is defined as.\\ $\underline{\zeta_{\widetilde{p}}}(z_{1})=\bigwedge\limits_{z_{1}^{'}\in[z_{1}]_{\omega}}\zeta_{\widetilde{p}}(z_{1}^{'})$ and $\underline{\eta_{\widetilde{p}}}(z_{1})=\bigvee\limits_{z_{1}^{'}\in[z_{1}]_{\omega}}\eta_{\widetilde{p}}(z_{1}^{'})$, 
$\underline{\zeta_{p}}(z_{1})=\bigwedge\limits_{z_{1}^{'}\in[z_{1}]_{\omega}}\zeta_{p}(z_{1}^{'})$ and $\underline{\eta_{p}}(z_{1})=\bigvee\limits_{z_{1}^{'}\in[z_{1}]_{\omega}}\eta_{p}(z_{1}^{'})$ Now, for all $z_{1}\in[z_{1}]_{\omega}$,\\ we have $\left(\zeta_{p^{\Box}}(z_{1})\right)^{2}+\left(\eta_{p^{\Box}}(z_{1})\right)^{2}\\=\left\{\left\langle\left[\bigwedge\limits_{z_{1}^{'}\in[z_{1}]_{\omega}}\zeta_{\widetilde{p}}(z_{1}^{'})\right]^{2}+\left[\bigvee\limits_{z_{1}^{'}\in[z_{1}]_{\omega}}\eta_{\widetilde{p}}(z_{1}^{'})\right]^{2}\right\rangle,\left(\bigwedge\limits_{z_{1}^{'}\in[z_{1}]_{\omega}}\zeta_{p}(z_{1}^{'})\right)^{2}+\left(\bigvee\limits_{z_{1}^{'}\in[z_{1}]_{\omega}}\eta_{p}(z_{1}^{'})\right)^{2}\right\}$\\
$=\left\langle\left[\bigwedge\limits_{z_{1}^{'}\in[z_{1}]_{\omega}}\zeta_{\widetilde{p}}(z_{1}^{'}),\bigwedge\limits_{z_{1}^{'}\in[z_{1}]_{\omega}}\zeta_{p}(z_{1}^{'})\right]^{2}+\left[\bigvee\limits_{z_{1}^{'}\in[z_{1}]_{\omega}}\eta_{\widetilde{p}}(z_{1}^{'}),\bigvee\limits_{z_{1}^{'}\in[z_{1}]_{\omega}}\eta_{p}(z_{1}^{'})\right]^{2}\right\rangle$\\
$=\bigwedge\limits_{z_{1}^{'}\in[z_{1}]_{\omega}}\left[\zeta_{p^{\Box}}(z_{1}^{'})\right]^{2}+\bigvee\limits_{z_{1}^{'}\in[z_{1}]_{\omega}}\left[\eta_{p^{\Box}}(z_{1}^{'})\right]^{2}$\\
$\leq\bigwedge\limits_{z_{1}^{'}\in[z_{1}]_{\omega}}\left(\zeta_{p^{\Box}}(z_{1}^{'})\right)^{2}+\bigvee\limits_{z_{1}^{'}\in[z_{1}]_{\omega}}\left(1-\left(\eta_{p^{\Box}}(z_{1}^{'})\right)^{2}\right)$\\
$=\bigwedge\limits_{z_{1}^{'}\in[z_{1}]_{\omega}}\left(\zeta_{p^{\Box}}(z_{1}^{'})\right)^{2}+1-\bigwedge\limits_{z_{1}^{'}\in[z_{1}]_{\omega}}\left(\eta_{p^{\Box}}(z_{1}^{'})\right)^{2}$\\
implies $\left(\zeta_{p^{\Box}}(z_{1})\right)^{2}+\left(\eta_{p^{\Box}}(z_{1})\right)^{2}\leq 1$\\
Similarly, $\overline{App}(P^{\Box})$.
\end{proof}
\end{proposition}
\begin{theorem}
Let us consider any two CPFSs $P^{\Box}_{1}=\left\langle\left[\zeta_{\widetilde{p_{1}}},\eta_{\widetilde{p_{1}}}\right],\left(\zeta_{p_{1}},\eta_{p_{1}}\right)\right\rangle$ and $P^{\Box}_{2}=\left\langle\left[\zeta_{\widetilde{p_{2}}},\eta_{\widetilde{p_{2}}}\right],\left(\zeta_{p_{2}},\eta_{p_{2}}\right)\right\rangle$ of $S$ and $\omega$ be the complete $CR_{\omega}$ on $S$. Then $\underline{App}(P^{\Box}_{1})\circ\underline{App}(P^{\Box}_{2})\subseteq \underline{App}\left(P^{\Box}_{1}\circ P^{\Box}_{2}\right)$
\begin{proof}
Since $\omega$ is a complete $CR_{\omega}$ on $S$ so $[z_{1}]_{\omega}[z_{2}]_{\omega}=[z_{1}z_{2}]_{\omega}$ for all $z_{1},z_{2}\in S$ As $\underline{App}(P^{\Box}_{1})=\left\langle\left[\zeta_{\widetilde{p_{1}}},\eta_{\widetilde{p_{1}}}\right],\left(\zeta_{p_{1}},\eta_{p_{1}}\right)\right\rangle$ and $\underline{App}(P^{\Box}_{2})=\left\langle\left[\zeta_{\widetilde{p_{2}}},\eta_{\widetilde{p_{2}}}\right],\left(\zeta_{p_{2}},\eta_{p_{2}}\right)\right\rangle$. Now, $\underline{App}(P_{1}^{\Box})\circ\underline{App}(P_{2}^{\Box})=\left(\zeta_{p^{\Box}_{1}}\circ\zeta_{p^{\Box}_{2}},\eta_{p^{\Box}_{1}}\circ\eta_{p^{\Box}_{2}}\right)$ and $\underline{App}\left(P^{\Box}_{1}\circ P^{\Box}_{2}\right)=\left(\left(\zeta_{p^{\Box}_{1}}\circ\zeta_{p^{\Box}_{2}}\right),\left(\eta_{p^{\Box}_{1}}\circ\eta_{p^{\Box}_{2}}\right)\right)$.To show that $\underline{App}(P^{\Box}_{1})\circ\underline{App}(P^{\Box}_{2})\subseteq \underline{App}\left(P^{\Box}_{1}\circ P^{\Box}_{2}\right)$, we have to prove that $\left[\underline{\zeta_{p^{\Box}_{1}}}\circ\underline{\zeta_{p^{\Box}_{2}}}\right](z_{1})\leq\left(\underline{\zeta_{p^{\Box}_{1}}\circ\zeta_{p^{\Box}_{2}}}\right)(z_{1})$ and $\left[\underline{\eta_{p^{\Box}_{1}}}\circ\underline{\eta_{p^{\Box}_{2}}}\right](z_{1})\geq\left(\underline{\eta_{p^{\Box}_{1}}\circ\eta_{p^{\Box}_{2}}}\right)(z_{1})$ Now, for all $z\in S$\\ $\left[\underline{\zeta_{p^{\Box}_{1}}}\circ\underline{\zeta_{p^{\Box}_{2}}}\right](z)=\bigvee\limits_{z=z_{1}z_{2}}\left(\underline{\zeta_{p^{\Box}_{1}}(z_{1})}\bigwedge\underline{\zeta_{p^{\Box}_{2}}(z_{2})}\right)$\\
$\mbox{\hspace{1.3cm}}=\bigvee\limits_{z=z_{1}z_{2}}\left[\left(\underline{\zeta_{\widetilde{p}_{1}}}(z_{1}),\underline{\zeta_{p_{1}}}(z_{1})\right)\bigwedge\left(\underline{\zeta_{\widetilde{p}_{2}}}(z_{2}),\underline{\zeta_{p_{2}}}(z_{2})\right)\right]$\\
$\mbox{\hspace{1.3cm}}=\bigvee\limits_{z=z_{1}z_{2}}\left[\left(\bigwedge\limits_{M\in[z_{1}]_{\omega}}\zeta_{\widetilde{p_{1}}}(M),\bigwedge\limits_{M\in[z_{1}]_{\omega}}\zeta_{p_{1}}(M)\right)\bigwedge\left(\bigwedge\limits_{N\in[z_{2}]_{\omega}}\zeta_{\widetilde{p_{2}}}(N),\bigwedge\limits_{N\in[z_{2}]_{\omega}}\zeta_{p_{2}}(N)\right)\right]$\\
$\mbox{\hspace{1.3cm}}=\bigvee\limits_{z=z_{1}z_{2}}\left[\left(\bigwedge\limits_{M\in[z_{1}]_{\omega}}\zeta_{p_{1}^{\Box}}(M)\right)\bigwedge\left(\bigwedge\limits_{N\in[z_{2}]_{\omega}}\zeta_{p_{2}^{\Box}}(N)\right)\right]$\\
$\mbox{\hspace{1.3cm}}=\bigvee\limits_{z=z_{1}z_{2}}\left[\bigwedge\limits_{M\in[z_{1}]_{\omega}N\in[z_{2}]_{\omega}}\left(\zeta_{p_{1}^{\Box}}(M)\bigwedge\zeta_{p_{2}^{\Box}}(N)\right)\right]$\\
$\mbox{\hspace{1.3cm}}\leq\bigvee\limits_{z=z_{1}z_{2}}\left[\bigwedge\limits_{MN\in[z_{1}z_{2}]_{\omega}}\left(\zeta_{p_{1}^{\Box}}(M)\bigwedge\zeta_{p_{2}^{\Box}}(N)\right)\right]$ $\mbox{\hspace{1.2cm}}$as $MN\in[z_{1}]_{\omega}[z_{2}]_{\omega}=[z_{1}z_{2}]_{\omega}$
$\mbox{\hspace{1.3cm}}=\bigvee\limits_{MN\in[z]_{\omega}}\left(\zeta_{p_{1}^{\Box}}(M)\bigwedge\zeta_{p_{2}^{\Box}}(N)\right)$\\
$\mbox{\hspace{1.3cm}}=\bigvee\limits_{\lambda\in[z]_{\omega},\lambda=MN}\left(\zeta_{p_{1}^{\Box}}(M)\bigwedge\zeta_{p_{2}^{\Box}}(N)\right)$\\
$\mbox{\hspace{1.3cm}}=\bigvee\limits_{\lambda\in[z]_{\omega}}\left[\bigvee\limits_{\lambda=MN}\left(\zeta_{p_{1}^{\Box}}(M)\bigwedge\zeta_{p_{2}^{\Box}}(N)\right)\right]$\\
$\mbox{\hspace{1.3cm}}=\bigvee\limits_{\lambda\in[z]_{\omega}}\left[\left(\zeta_{p_{1}^{\Box}}\circ\zeta_{p_{2}^{\Box}}\right)(\lambda)\right]$\\
implies $\left[\underline{\zeta_{p_{1}^{\Box}}}\circ\underline{\zeta_{p_{2}^{\Box}}}\right](z)\leq\left[\underline{\zeta_{p_{1}^{\Box}}\circ\zeta_{p_{2}^{\Box}}}\right](z)$.\\
Further\\
$\left[\underline{\eta_{p^{\Box}_{1}}}\circ\underline{\eta_{p^{\Box}_{2}}}\right](z)=\bigwedge\limits_{z=z_{1}z_{2}}\left(\underline{\eta_{p^{\Box}_{1}}(z_{1})}\bigvee\underline{\eta_{p^{\Box}_{2}}(z_{2})}\right)$\\
$\mbox{\hspace{1.3cm}}=\bigwedge\limits_{z=z_{1}z_{2}}\left[\left(\underline{\eta_{\widetilde{p}_{1}}}(z_{1}),\underline{\eta_{p_{1}}}(z_{1})\right)\bigvee\left(\underline{\eta_{\widetilde{p}_{2}}}(z_{2}),\underline{\eta_{p_{2}}}(z_{2})\right)\right]$\\
$\mbox{\hspace{1.3cm}}=\bigwedge\limits_{z=z_{1}z_{2}}\left[\left(\bigvee\limits_{M\in[z_{1}]_{\omega}}\eta_{\widetilde{p_{1}}}(M),\bigvee\limits_{M\in[z_{1}]_{\omega}}\eta_{p_{1}}(M)\right)\bigvee\left(\bigvee\limits_{N\in[z_{2}]_{\omega}}\eta_{\widetilde{p_{2}}}(N),\bigvee\limits_{N\in[z_{2}]_{\omega}}\eta_{p_{2}}(N)\right)\right]$\\
$\mbox{\hspace{1.3cm}}=\bigwedge\limits_{z=z_{1}z_{2}}\left[\left(\bigvee\limits_{M\in[z_{1}]_{\omega}}\eta_{p_{1}^{\Box}}(M)\right)\bigvee\left(\bigvee\limits_{N\in[z_{2}]_{\omega}}\eta_{p_{2}^{\Box}}(N)\right)\right]$\\
$\mbox{\hspace{1.3cm}}=\bigwedge\limits_{z=z_{1}z_{2}}\left[\bigvee\limits_{M\in[z_{1}]_{\omega}N\in[z_{2}]_{\omega}}\left(\eta_{p_{1}^{\Box}}(M)\bigvee\eta_{p_{2}^{\Box}}(N)\right)\right]$\\
$\mbox{\hspace{1.3cm}}\geq\bigwedge\limits_{z=z_{1}z_{2}}\left[\bigvee\limits_{MN\in[z_{1}z_{2}]_{\omega}}\left(\eta_{p_{1}^{\Box}}(M)\bigvee\eta_{p_{2}^{\Box}}(N)\right)\right]$ $\mbox{\hspace{1.2cm}}$as $MN\in[z_{1}]_{\omega}[z_{2}]_{\omega}=[z_{1}z_{2}]_{\omega}$
$\mbox{\hspace{1.3cm}}=\bigwedge\limits_{MN\in[z]_{\omega}}\left(\eta_{p_{1}^{\Box}}(M)\bigvee\eta_{p_{2}^{\Box}}(N)\right)$\\
$\mbox{\hspace{1.3cm}}=\bigwedge\limits_{\lambda\in[z]_{\omega},\lambda=MN}\left(\eta_{p_{1}^{\Box}}(M)\bigvee\eta_{p_{2}^{\Box}}(N)\right)$\\
$\mbox{\hspace{1.3cm}}=\bigwedge\limits_{\lambda\in[z]_{\omega}}\left[\bigwedge\limits_{\lambda=MN}\left(\eta_{p_{1}^{\Box}}(M)\bigvee\eta_{p_{2}^{\Box}}(N)\right)\right]$\\
$\mbox{\hspace{1.3cm}}=\bigwedge\limits_{\lambda\in[z]_{\omega}}\left[\left(\eta_{p_{1}^{\Box}}\circ\eta_{p_{2}^{\Box}}\right)(\lambda)\right]$\\
implies $\left[\underline{\eta_{p_{1}^{\Box}}}\circ\underline{\eta_{p_{2}^{\Box}}}\right](z)\geq\left[\underline{\eta_{p_{1}^{\Box}}\circ\eta_{p_{2}^{\Box}}}\right](z)$.\\
Hence, $\underline{App}(P_{1}^{\Box})\circ\underline{App}(P_{2}^{\Box})\subseteq\underline{App}(P_{1}^{\Box}\circ P_{2}^{\Box})$.
\end{proof}
\end{theorem}
\begin{theorem}
Let $P_{1}^{\Box}$,$P_{2}^{\Box}$ be any two CPFSs on $S$. Then $P^{\Box}_{1}=\left\langle\left[\zeta_{\widetilde{p_{1}}},\eta_{\widetilde{p_{1}}}\right],\left(\zeta_{p_{1}},\eta_{p_{1}}\right)\right\rangle$ and $P^{\Box}_{2}=\left\langle\left[\zeta_{\widetilde{p_{2}}},\eta_{\widetilde{p_{2}}}\right],\left(\zeta_{p_{2}},\eta_{p_{2}}\right)\right\rangle$ of $S$ and let $\omega$ be the complete $CR_{\omega}$ on $S$. Then $\overline{App}(P^{\Box}_{1})\circ\overline{App}(P^{\Box}_{2})\subseteq \overline{App}\left(P^{\Box}_{1}\circ P^{\Box}_{2}\right)$
\begin{proof}
Since $\omega$ is a complete $CR_{\omega}$ on $S$ so $[z_{1}]_{\omega}[z_{2}]_{\omega}\subseteq[z_{1}z_{2}]_{\omega}$ for all $z_{1},z_{2}\in S$ As $\overline{App}(P^{\Box}_{1})=\left\langle\left[\zeta_{\widetilde{p_{1}}},\eta_{\widetilde{p_{1}}}\right],\left(\zeta_{p_{1}},\eta_{p_{1}}\right)\right\rangle$ and $\overline{App}(P^{\Box}_{2})=\left\langle\left[\zeta_{\widetilde{p_{2}}},\eta_{\widetilde{p_{2}}}\right],\left(\zeta_{p_{2}},\eta_{p_{2}}\right)\right\rangle$. Now, $\overline{App}(P_{1}^{\Box})\circ\overline{App}(P_{2}^{\Box})=\left(\zeta_{p^{\Box}_{1}}\circ\zeta_{p^{\Box}_{2}},\eta_{p^{\Box}_{1}}\circ\eta_{p^{\Box}_{2}}\right)$ and $\overline{App}\left(P^{\Box}_{1}\circ P^{\Box}_{2}\right)=\left(\left(\zeta_{p^{\Box}_{1}}\circ\zeta_{p^{\Box}_{2}}\right),\left(\eta_{p^{\Box}_{1}}\circ\eta_{p^{\Box}_{2}}\right)\right)$.To show that $\overline{App}(P^{\Box}_{1})\circ\overline{App}(P^{\Box}_{2})\subseteq \overline{App}\left(P^{\Box}_{1}\circ P^{\Box}_{2}\right)$, we have to prove that $\left[\overline{\zeta_{p^{\Box}_{1}}}\circ\overline{\zeta_{p^{\Box}_{2}}}\right](z_{1})\leq\left(\overline{\zeta_{p^{\Box}_{1}}\circ\zeta_{p^{\Box}_{2}}}\right)(z_{1})$ and $\left[\overline{\eta_{p^{\Box}_{1}}}\circ\overline{\eta_{p^{\Box}_{2}}}\right](z_{1})\geq\left(\overline{\eta_{p^{\Box}_{1}}\circ\eta_{p^{\Box}_{2}}}\right)(z_{1})$ Now, for all $z\in S$\\ $\left[\overline{\zeta_{p^{\Box}_{1}}}\circ\overline{\zeta_{p^{\Box}_{2}}}\right](z)=\bigvee\limits_{z=z_{1}z_{2}}\left(\overline{\zeta_{p^{\Box}_{1}}}(z_{1})\bigwedge\overline{\zeta_{p^{\Box}_{2}}}(z_{2})\right)$\\
$\mbox{\hspace{1.3cm}}=\bigvee\limits_{z=z_{1}z_{2}}\left[\left(\overline{\zeta_{\widetilde{p}_{1}}}(z_{1}),\overline{\zeta_{p_{1}}}(z_{1})\right)\bigwedge\left(\overline{\zeta_{\widetilde{p}_{2}}}(z_{2}),\overline{\zeta_{p_{2}}}(z_{2})\right)\right]$\\
$\mbox{\hspace{1.3cm}}=\bigvee\limits_{z=z_{1}z_{2}}\left[\left(\bigvee\limits_{M\in[z_{1}]_{\omega}}\zeta_{\widetilde{p_{1}}}(M),\bigvee\limits_{M\in[z_{1}]_{\omega}}\zeta_{p_{1}}(M)\right)\bigwedge\left(\bigvee\limits_{N\in[z_{2}]_{\omega}}\zeta_{\widetilde{p_{2}}}(N),\bigvee\limits_{N\in[z_{2}]_{\omega}}\zeta_{p_{2}}(N)\right)\right]$\\
$\mbox{\hspace{1.3cm}}=\bigvee\limits_{z=z_{1}z_{2}}\left[\left(\bigvee\limits_{M\in[z_{1}]_{\omega}}\zeta_{p_{1}^{\Box}}(M)\right)\bigwedge\left(\bigvee\limits_{N\in[z_{2}]_{\omega}}\zeta_{p_{2}^{\Box}}(N)\right)\right]$\\
$\mbox{\hspace{1.3cm}}=\bigvee\limits_{z=z_{1}z_{2}}\left[\bigvee\limits_{M\in[z_{1}]_{\omega}N\in[z_{2}]_{\omega}}\left(\zeta_{p_{1}^{\Box}}(M)\bigwedge\zeta_{p_{2}^{\Box}}(N)\right)\right]$\\
$\mbox{\hspace{1.3cm}}\leq\bigvee\limits_{z=z_{1}z_{2}}\left[\bigvee\limits_{MN\in[z_{1}z_{2}]_{\omega}}\left(\zeta_{p_{1}^{\Box}}(M)\bigwedge\zeta_{p_{2}^{\Box}}(N)\right)\right]$ $\mbox{\hspace{1.2cm}}$as $MN\in[z_{1}]_{\omega}[z_{2}]_{\omega}=[z_{1}z_{2}]_{\omega}$\\
$\mbox{\hspace{1.3cm}}=\bigvee\limits_{MN\in[z]_{\omega}}\left(\zeta_{p_{1}^{\Box}}(M)\bigwedge\zeta_{p_{2}^{\Box}}(N)\right)$\\
$\mbox{\hspace{1.3cm}}=\bigvee\limits_{\lambda\in[z]_{\omega},\lambda=MN}\left(\zeta_{p_{1}^{\Box}}(M)\bigwedge\zeta_{p_{2}^{\Box}}(N)\right)$\\
$\mbox{\hspace{1.3cm}}=\bigvee\limits_{\lambda\in[z]_{\omega}}\left[\bigvee\limits_{\lambda=MN}\left(\zeta_{p_{1}^{\Box}}(M)\bigwedge\zeta_{p_{2}^{\Box}}(N)\right)\right]$\\
$\mbox{\hspace{1.3cm}}=\bigvee\limits_{\lambda\in[z]_{\omega}}\left[\left(\zeta_{p_{1}^{\Box}}\circ\zeta_{p_{2}^{\Box}}\right)(\lambda)\right]$\\
implies $\left[\overline{\zeta_{p_{1}^{\Box}}}\circ\overline{\zeta_{p_{2}^{\Box}}}\right](z)\leq\left[\overline{\zeta_{p_{1}^{\Box}}\circ\zeta_{p_{2}^{\Box}}}\right](z)$.\\
Further\\
$\left[\overline{\eta_{p^{\Box}_{1}}}\circ\overline{\eta_{p^{\Box}_{2}}}\right](z)=\bigwedge\limits_{z=z_{1}z_{2}}\left(\overline{\eta_{p^{\Box}_{1}}(z_{1})}\bigvee\overline{\eta_{p^{\Box}_{2}}(z_{2})}\right)$\\
$\mbox{\hspace{1.3cm}}=\bigwedge\limits_{z=z_{1}z_{2}}\left[\left(\overline{\eta_{\widetilde{p}_{1}}}(z_{1}),\overline{\eta_{p_{1}}}(z_{1})\right)\bigvee\left(\overline{\eta_{\widetilde{p}_{2}}}(z_{2}),\overline{\eta_{p_{2}}}(z_{2})\right)\right]$\\
$\mbox{\hspace{1.3cm}}=\bigwedge\limits_{z=z_{1}z_{2}}\left[\left(\bigwedge\limits_{M\in[z_{1}]_{\omega}}\eta_{\widetilde{p_{1}}}(M),\bigwedge\limits_{M\in[z_{1}]_{\omega}}\eta_{p_{1}}(M)\right)\bigvee\left(\bigwedge\limits_{N\in[z_{2}]_{\omega}}\eta_{\widetilde{p_{2}}}(N),\bigwedge\limits_{N\in[z_{2}]_{\omega}}\eta_{p_{2}}(N)\right)\right]$\\
$\mbox{\hspace{1.3cm}}=\bigwedge\limits_{z=z_{1}z_{2}}\left[\left(\bigwedge\limits_{M\in[z_{1}]_{\omega}}\eta_{p_{1}^{\Box}}(M)\right)\bigvee\left(\bigwedge\limits_{N\in[z_{2}]_{\omega}}\eta_{p_{2}^{\Box}}(N)\right)\right]$\\
$\mbox{\hspace{1.3cm}}=\bigwedge\limits_{z=z_{1}z_{2}}\left[\bigwedge\limits_{M\in[z_{1}]_{\omega}N\in[z_{2}]_{\omega}}\left(\eta_{p_{1}^{\Box}}(M)\bigvee\eta_{p_{2}^{\Box}}(N)\right)\right]$\\
$\mbox{\hspace{1.3cm}}\geq\bigwedge\limits_{z=z_{1}z_{2}}\left[\bigwedge\limits_{MN\in[z_{1}z_{2}]_{\omega}}\left(\eta_{p_{1}^{\Box}}(M)\bigvee\eta_{p_{2}^{\Box}}(N)\right)\right]$ $\mbox{\hspace{1.2cm}}$as $MN\in[z_{1}]_{\omega}[z_{2}]_{\omega}=[z_{1}z_{2}]_{\omega}$\\
$\mbox{\hspace{1.3cm}}=\bigwedge\limits_{MN\in[z]_{\omega}}\left(\eta_{p_{1}^{\Box}}(M)\bigvee\eta_{p_{2}^{\Box}}(N)\right)$\\
$\mbox{\hspace{1.3cm}}=\bigwedge\limits_{\lambda\in[z]_{\omega},\lambda=MN}\left(\eta_{p_{1}^{\Box}}(M)\bigvee\eta_{p_{2}^{\Box}}(N)\right)$\\
$\mbox{\hspace{1.3cm}}=\bigwedge\limits_{\lambda\in[z]_{\omega}}\left[\bigwedge\limits_{\lambda=MN}\left(\eta_{p_{1}^{\Box}}(M)\bigvee\eta_{p_{2}^{\Box}}(N)\right)\right]$\\
$\mbox{\hspace{1.3cm}}=\bigwedge\limits_{\lambda\in[z]_{\omega}}\left[\left(\eta_{p_{1}^{\Box}}\circ\eta_{p_{2}^{\Box}}\right)(\lambda)\right]$\\
implies $\left[\overline{\eta_{p_{1}^{\Box}}}\circ\overline{\eta_{p_{2}^{\Box}}}\right](z)\geq\left[\overline{\eta_{p_{1}^{\Box}}\circ\eta_{p_{2}^{\Box}}}\right](z)$.\\
Hence, $\overline{App}(P_{1}^{\Box})\circ\overline{App}(P_{2}^{\Box})\subseteq\overline{App}(P_{1}^{\Box}\circ P_{2}^{\Box})$.
\end{proof}
\end{theorem}
\section{Rough cubic Pythagorean fuzzy ideals (RCPFI) in semigroup.}
In this section, $P^{\Box}_{LI}$, $P^{\Box}_{RI}$, $P^{\Box}_{I}$, $P^{\Box}_{BI}$, $P^{\Box}_{II}$ cubic Pythagorean fuzzy left ideal, cubic Pythagorean fuzzy right ideal, cubic Pythagorean fuzzy ideal, cubic Pythagorean fuzzy bi-ideal and cubic Pythagorean fuzzy interior-ideal are respectively.

\begin{definition}
Let $\omega$ be a $CR_{\omega}$ on $S$ and $P^{\Box}$ be a CPFS. Then $P^{\Box}$ is called lower (resp.upper) rough cubic Pythagorean fuzzy sub-semigroup of $S$, if $\underline{App}(P^{\Box})$ (resp.$\overline{App}(P^{\Box})$) is a cubic Pythagorean fuzzy sub-semigroup of $S$.\\ A cubic Pythagorean fuzzy set $P^{\Box}$ is known to be rough cubic Pythagorean fuzzy sub-semigroup of $S$, if $\underline{App}(P^{\Box})$ and $\overline{App}(P^{\Box})$ are both Pythagorean fuzzy sub-semigroup of $S$.
\end{definition}
\begin{definition}
Let $\omega$ be a $CR_{\omega}$ on $S$ and $P^{\Box}$ be a CPFS. Then $P^{\Box}$ is called lower rough $P_{LI}^{\Box}$ (resp.$P_{RI}^{\Box}$,$P_{I}^{\Box}$) of $S$, if $\underline{App}(P^{\Box})$ is a $P_{LI}^{\Box}$ (resp.$P_{RI}^{\Box}$,$P_{I}^{\Box}$) of $S$ and\\$(i)~~\underline{\zeta_{\widetilde{p}}}(xy)\geq\underline{\zeta_{\widetilde{p}}}(y)$; $\underline{\zeta_{p}}(xy)\leq\underline{\zeta_{p}}(y)~~~\forall~x,y\in~S$\\$(ii)~~\underline{\eta_{\widetilde{p}}}(xy)\geq\underline{\eta_{\widetilde{p}}}(y)$;$\underline{\eta_{p}}(xy)\leq\underline{\eta_{p}}(y)~~~\forall~x,y\in~S$
\end{definition}
\begin{definition}
Let $\omega$ be a $CR_{\omega}$ on $S$ and $P^{\Box}$ be a CPFS. Then $P^{\Box}$ is called upper rough $P_{LI}^{\Box}$ (resp.$P_{RI}^{\Box}$,$P_{I}^{\Box}$) of $S$, if $\overline{App}(P^{\Box})$ is a $P_{LI}^{\Box}$ (resp.$P_{RI}^{\Box}$,$P_{I}^{\Box}$) of $S$ and\\$(i)~~\overline{\zeta_{\widetilde{p}}}(xy)\geq\overline{\zeta_{\widetilde{p}}}(y)$; $\overline{\zeta_{p}}(xy)\leq\overline{\zeta_{p}}(y)~~~\forall~x,y\in~S$\\$(ii)~~\overline{\eta_{\widetilde{p}}}(xy)\geq\overline{\eta_{\widetilde{p}}}(y)$;$\overline{\eta_{p}}(xy)\leq\overline{\eta_{p}}(y)~~~\forall~x,y\in~S$
\end{definition}
\begin{definition}
Let $P^{\Box}$ be a CPFS and $\omega$ be a $CR_{\omega}$ on $S$. Then $P^{\Box}$ is called lower(resp. upper) rough $P_{BI}^{\Box}$ of $S$, if $\underline{App}(P^{\Box})$ (resp. $\overline{App}(P^{\Box})$) is a $P_{BI}^{\Box}$ of $S$ and \\
$(i)~~~~~\zeta_{\widetilde{p}}(xyz)\geq min\left\{\zeta_{\widetilde{p}}(x),\zeta_{\widetilde{p}}(z)\right\}~~~\forall x,y,z\in~S.$\\
$(ii)~~~~\eta_{\widetilde{p}}(xyz)\geq min\left\{\eta_{\widetilde{p}}(x),\eta_{\widetilde{p}}(z)\right\}~~~\forall x,y,z\in~S.$\\
$(iii)~~~\zeta_{p}(xyz)\leq max\left\{\zeta_{p}(x),\zeta_{p}(z)\right\}~~~\forall x,y,z\in~S.$\\
$(iv)~\eta_{p}(xyz)\leq max\left\{\eta_{p}(x),\eta_{p}(z)\right\}~~~\forall x,y,z\in~S.$
\end{definition}
\begin{definition}
Let $P^{\Box}$ be a CPFS and $\omega$ be a $CR_{\omega}$ on $S$. Then $P^{\Box}$ is called lower(resp. upper) rough $P_{II}^{\Box}$ of $S$, if $\underline{App}(P^{\Box})$ (resp. $\overline{App}(P^{\Box})$) is a $P_{II}^{\Box}$ of $S$ and \\
$(i)~~~~~\zeta_{\widetilde{p}}(xyz)\geq \zeta_{\widetilde{p}}(y)~~~\forall x,y,z\in~S.$\\
$(ii)~~~~\eta_{\widetilde{p}}(xyz)\geq \eta_{\widetilde{p}}(y)~~~\forall x,y,z\in~S.$\\
$(iii)~~~\zeta_{p}(xyz)\leq \zeta_{p}(y)~~~\forall x,y,z\in~S.$\\
$(iv)~\eta_{p}(xyz)\leq \eta_{p}(y)~~~\forall x,y,z\in~S.$
\end{definition}
\begin{theorem}
Let $\omega$ is a $CR_{\omega}$ on $S$ and $P^{\Box}$ be a cubic Pythagorean fuzzy sub-semigroup of $S$. Then $\overline{App}(P^{\Box})$ is a cubic Pythagorean fuzzy sub-semigroup of $S$.
\begin{proof}
Since $\omega$ is a $CR_{\omega}$ on $S$, for all $z_{1},z_{2}\in S$, we have $[z_{1}][z_{2}]\subseteq[z_{1}z_{2}]_{\omega}$. Now, we have to show that $\overline{App}(P^{\Box})=\left(\overline{\zeta_{p^{\Box}}},\overline{\eta_{p^{\Box}}}\right)$ is a cubic Pythagorean fuzzy sub-semigroup of $S$, consider\\
$\overline{\zeta_{\widetilde{p}}}(z_{1},z_{2})=\bigvee\limits_{z_{3}\in[z_{1}z_{2}]_{\omega}}\zeta_{\widetilde{p}}(z_{3})\geq\bigvee\limits_{z_{3}\in[z_{1}]_{\omega}[z_{2}]_{\omega}}\zeta_{\widetilde{p}}(z_{3})$\\$\mbox{\hspace{1.3cm}}=\bigvee\limits_{MN\in[z_{1}]_{\omega}[z_{2}]_{\omega}}\zeta_{\widetilde{p}}(MN)$\\$\mbox{\hspace{1.3cm}}\geq\bigvee\limits_{M\in[z_{1}]_{\omega},N\in[z_{2}]_{\omega}}\left[\zeta_{\widetilde{p}}(M)\bigwedge\zeta_{\widetilde{p}}(N)\right]$\\
$\mbox{\hspace{1.3cm}}=\left[\bigvee\limits_{M\in[z_{1}]_{\omega}}\left[\zeta_{\widetilde{p}}(M)\right]\right]\bigwedge\left[\bigvee\limits_{N\in[z_{2}]_{\omega}}\zeta_{\widetilde{p}}(N)\right]$\\
implies $\overline{\zeta_{\widetilde{p}}}(z_{1},z_{2})\geq min\left\{\overline{\zeta_{\widetilde{p}}}(z_{1}),\overline{\zeta_{\widetilde{p}}}(z_{2})\right\}$\\
$\overline{\zeta_{p}}(z_{1},z_{2})=\bigvee\limits_{z_{3}\in[z_{1}]_{\omega}[z_{2}]_{\omega}}\zeta_{p}(z_{3})$\\
$\mbox{\hspace{1.3cm}}\leq\bigvee\limits_{z_{3}\in[z_{1}z_{2}]_{\omega}}\zeta_{p}(z_{3})$\\$\mbox{\hspace{1.3cm}}=\bigvee\limits_{MN\in[z_{1}]_{\omega}[z_{2}]_{\omega}}\zeta_{p}(MN)$\\$\mbox{\hspace{1.3cm}}\leq\bigvee\limits_{MN\in[z_{1}z_{2}]_{\omega}}\zeta_{p}(MN)$\\
$\mbox{\hspace{1.3cm}}=\left(\bigvee\limits_{M\in[z_{1}]_{\omega}}\zeta_{p}(M)\right)\bigvee\left(\bigvee\limits_{N\in[z_{2}]_{\omega}}\zeta_{p}(N)\right)$\\$\overline{\zeta_{p}}(z_{1}z_{2})\leq max\left\{\overline{\zeta_{p}}(M),\overline{\zeta_{p}}(N)\right\}$\\
Further\\
$\overline{\eta_{\widetilde{p}}}(z_{1},z_{2})=\bigwedge\limits_{z_{3}\in[z_{1}z_{2}]_{\omega}}\eta_{\widetilde{p}}(z_{3})\geq\bigwedge\limits_{z_{3}\in[z_{1}]_{\omega}[z_{2}]_{\omega}}\eta_{\widetilde{p}}(z_{3})$\\$\mbox{\hspace{1.3cm}}=\bigwedge\limits_{MN\in[z_{1}]_{\omega}[z_{2}]_{\omega}}\eta_{\widetilde{p}}(MN)$\\$\mbox{\hspace{1.3cm}}\geq\bigwedge\limits_{M\in[z_{1}]_{\omega},N\in[z_{2}]_{\omega}}\left[\eta_{\widetilde{p}}(M)\bigwedge\eta_{\widetilde{p}}(N)\right]$\\
$\mbox{\hspace{1.3cm}}=\left[\bigwedge\limits_{M\in[z_{1}]_{\omega}}\left[\eta_{\widetilde{p}}(M)\right]\right]\bigwedge\left[\bigwedge\limits_{N\in[z_{2}]_{\omega}}\eta_{\widetilde{p}}(N)\right]$\\
implies $\overline{\eta_{\widetilde{p}}}(z_{1},z_{2})\geq min\left\{\overline{\eta_{\widetilde{p}}}(z_{1}),\overline{\eta_{\widetilde{p}}}(z_{2})\right\}$\\
$\overline{\eta_{p}}(z_{1},z_{2})=\bigwedge\limits_{z_{3}\in[z_{1}]_{\omega}[z_{2}]_{\omega}}\eta_{p}(z_{3})$\\
$\mbox{\hspace{1.3cm}}\leq\bigwedge\limits_{z_{3}\in[z_{1}z_{2}]_{\omega}}\eta_{p}(z_{3})$\\$\mbox{\hspace{1.3cm}}=\bigwedge\limits_{MN\in[z_{1}]_{\omega}[z_{2}]_{\omega}}\eta_{p}(MN)$\\$\mbox{\hspace{1.3cm}}\leq\bigwedge\limits_{MN\in[z_{1}z_{2}]_{\omega}}\eta_{p}(MN)$\\
$\mbox{\hspace{1.3cm}}=\left(\bigwedge\limits_{M\in[z_{1}]_{\omega}}\eta_{p}(M)\right)\bigvee\left(\bigwedge\limits_{N\in[z_{2}]_{\omega}}\eta_{p}(N)\right)$\\$\overline{\eta_{p}}(z_{1}z_{2})\leq max\left\{\overline{\eta_{p}}(M),\overline{\eta_{p}}(N)\right\}$.
\end{proof}
\end{theorem}
\begin{theorem}
Let $\omega$ be a $CR_{\omega}$ on $S$, and $P^{\Box}$ be a $P_{LI}^{\Box}$ (resp. $P_{RI}^{\Box}$) of $S$. Then $\overline{App}(P^{\Box})$ is a $P_{LI}^{\Box}$ (resp. $P_{RI}^{\Box}$) of $S$.
\begin{proof}
Since $\omega$ is a $CR_{\omega}$ on $S$, then for all $z_{1},z_{2}\in S$ it follows that $[z_{1}][z_{2}]\subseteq[z_{1}z_{2}]_{\omega}$. Now we have to show that $\overline{App}(P^{\Box})=\left(\overline{\zeta_{p^{\Box}}},\overline{\eta_{p^{\Box}}}\right)=\left\langle\left[\overline{\zeta_{\widetilde{p}}},\overline{\eta_{\widetilde{p}}}\right],\left(\overline{\zeta_{p}},\overline{\eta_{p}}\right)\right\rangle$ is a $P_{LI}^{\Box}$ of $S$.\\
$\overline{\zeta_{\widetilde{p}}}(z_{1}z_{2})=\bigvee\limits_{z_{3}\in[z_{1}z_{2}]_{\omega}}\zeta_{\widetilde{p}}(z_{3})\geq\bigvee\limits_{z_{3}\in[z_{1}]_{\omega}[z_{2}]_{\omega}}\zeta_{\widetilde{p}}(z_{3})$\\$\mbox{\hspace{1.3cm}}=\bigvee\limits_{MN\in[z_{1}]_{\omega}[z_{2}]_{\omega}}\zeta_{\widetilde{p}}(MN)\geq\bigvee\limits_{N\in[z_{2}]_{\omega}}\zeta_{\widetilde{p}}(N)$\\
implies $\overline{\zeta_{\widetilde{p}}}(z_{1}z_{2})\geq\overline{\zeta_{\widetilde{p}}}(z_{2})$\\
$\overline{\zeta_{p}}(z_{1}z_{2})=\bigvee\limits_{z_{3}\in[z_{1}z_{2}]_{\omega}}\zeta_{p}(z_{3})\leq\bigvee\limits_{z_{3}\in[z_{1}]_{\omega}[z_{2}]_{\omega}}\zeta_{p}(z_{3})$\\$\mbox{\hspace{1.3cm}}=\bigvee\limits_{MN\in[z_{1}]_{\omega}[z_{2}]_{\omega}}\zeta_{p}(MN)\leq\bigvee\limits_{N\in[z_{2}]_{\omega}}\zeta_{p}(N)$\\implies $\overline{\zeta_{p}}(z_{1}z_{2})\leq\overline{\zeta_{p}}(z_{2})$\\Next\\
$\overline{\eta_{\widetilde{p}}}(z_{1}z_{2})=\bigwedge\limits_{z_{3}\in[z_{1}z_{2}]_{\omega}}\eta_{\widetilde{p}}(z_{3})\geq\bigwedge\limits_{z_{3}\in[z_{1}]_{\omega}[z_{2}]_{\omega}}\eta_{\widetilde{p}}(z_{3})$\\$\mbox{\hspace{1.3cm}}=\bigwedge\limits_{MN\in[z_{1}]_{\omega}[z_{2}]_{\omega}}\eta_{\widetilde{p}}(MN)\geq\bigwedge\limits_{N\in[z_{2}]_{\omega}}\eta_{\widetilde{p}}(N)$\\
implies $\overline{\eta_{\widetilde{p}}}(z_{1}z_{2})\geq\overline{\eta_{\widetilde{p}}}(z_{2})$\\
$\overline{\eta_{p}}(z_{1}z_{2})=\bigwedge\limits_{z_{3}\in[z_{1}z_{2}]_{\omega}}\eta_{p}(z_{3})\leq\bigwedge\limits_{z_{3}\in[z_{1}]_{\omega}[z_{2}]_{\omega}}\eta_{p}(z_{3})$\\$\mbox{\hspace{1.3cm}}=\bigwedge\limits_{MN\in[z_{1}]_{\omega}[z_{2}]_{\omega}}\eta_{p}(MN)\leq\bigwedge\limits_{N\in[z_{2}]_{\omega}}\eta_{p}(N)$\\implies $\overline{\eta_{p}}(z_{1}z_{2})\leq\overline{\eta_{p}}(z_{2})$\\ implies that $\overline{App}(P^{\Box})$ is a $P_{LI}^{\Box}$ of $S$. Similarly, $\overline{App}(P^{\Box})$ is a $P_{RI}^{\Box}$ of $S$.
\end{proof}
\end{theorem}
\begin{theorem}
Let $\omega$ is a $CR_{\omega}$ on $S$ and let $P^{\Box}$ be a cubic Pythagorean fuzzy sub-semigroup of $S$. Then $\underline{App}(P^{\Box})$ is a cubic Pythagorean fuzzy sub-semigroup of $S$.
\begin{proof}
Since $\omega$ is a $CR_{\omega}$ on $S$, then for all $z_{1},z_{2}\in S$, $[z_{1}][z_{2}]=[z_{1}z_{2}]_{\omega}$. It is required to show that $\underline{App}(P^{\Box})=\left(\underline{\zeta_{p^{\Box}}},\underline{\eta_{p^{\Box}}}\right)$ is a cubic Pythagorean fuzzy sub-semigroup of $S$, consider\\
$\underline{\zeta_{\widetilde{p}}}(z_{1},z_{2})=\bigwedge\limits_{z_{3}\in[z_{1}z_{2}]_{\omega}}\zeta_{\widetilde{p}}(z_{3})\geq\bigwedge\limits_{z_{3}\in[z_{1}]_{\omega}[z_{2}]_{\omega}}\zeta_{\widetilde{p}}(z_{3})$\\$\mbox{\hspace{1.3cm}}=\bigwedge\limits_{MN\in[z_{1}]_{\omega}[z_{2}]_{\omega}}\zeta_{\widetilde{p}}(MN)$\\$\mbox{\hspace{1.3cm}}\geq\bigwedge\limits_{M\in[z_{1}]_{\omega},N\in[z_{2}]_{\omega}}\left[\zeta_{\widetilde{p}}(M)\bigwedge\zeta_{\widetilde{p}}(N)\right]$\\
$\mbox{\hspace{1.3cm}}=\left[\bigwedge\limits_{M\in[z_{1}]_{\omega}}\left[\zeta_{\widetilde{p}}(M)\right]\right]\bigwedge\left[\bigwedge\limits_{N\in[z_{2}]_{\omega}}\zeta_{\widetilde{p}}(N)\right]$\\
implies $\underline{\zeta_{\widetilde{p}}}(z_{1},z_{2})\geq min\left\{\underline{\zeta_{\widetilde{p}}}(z_{1}),\underline{\zeta_{\widetilde{p}}}(z_{2})\right\}$\\
$\underline{\zeta_{p}}(z_{1},z_{2})=\bigwedge\limits_{z_{3}\in[z_{1}]_{\omega}[z_{2}]_{\omega}}\zeta_{p}(z_{3})$\\
$\mbox{\hspace{1.3cm}}\leq\bigwedge\limits_{z_{3}\in[z_{1}z_{2}]_{\omega}}\zeta_{p}(z_{3})$\\$\mbox{\hspace{1.3cm}}=\bigwedge\limits_{MN\in[z_{1}]_{\omega}[z_{2}]_{\omega}}\zeta_{p}(MN)$\\$\mbox{\hspace{1.3cm}}\leq\bigwedge\limits_{MN\in[z_{1}z_{2}]_{\omega}}\zeta_{p}(MN)$\\
$\mbox{\hspace{1.3cm}}=\left(\bigwedge\limits_{M\in[z_{1}]_{\omega}}\zeta_{p}(M)\right)\bigwedge\left(\bigwedge\limits_{N\in[z_{2}]_{\omega}}\zeta_{p}(N)\right)$\\$\underline{\zeta_{p}}(z_{1}z_{2})\leq max\left\{\underline{\zeta_{p}}(M),\underline{\zeta_{p}}(N)\right\}$\\
Further\\
$\underline{\eta_{\widetilde{p}}}(z_{1},z_{2})=\bigvee\limits_{z_{3}\in[z_{1}z_{2}]_{\omega}}\eta_{\widetilde{p}}(z_{3})\geq\bigvee\limits_{z_{3}\in[z_{1}]_{\omega}[z_{2}]_{\omega}}\eta_{\widetilde{p}}(z_{3})$\\$\mbox{\hspace{1.3cm}}=\bigvee\limits_{MN\in[z_{1}]_{\omega}[z_{2}]_{\omega}}\eta_{\widetilde{p}}(MN)$\\$\mbox{\hspace{1.3cm}}\geq\bigvee\limits_{M\in[z_{1}]_{\omega},N\in[z_{2}]_{\omega}}\left[\eta_{\widetilde{p}}(M)\bigwedge\eta_{\widetilde{p}}(N)\right]$\\
$\mbox{\hspace{1.3cm}}=\left[\bigvee\limits_{M\in[z_{1}]_{\omega}}\left[\eta_{\widetilde{p}}(M)\right]\right]\bigwedge\left[\bigvee\limits_{N\in[z_{2}]_{\omega}}\eta_{\widetilde{p}}(N)\right]$\\
implies $\underline{\eta_{\widetilde{p}}}(z_{1},z_{2})\geq min\left\{\underline{\eta_{\widetilde{p}}}(z_{1}),\underline{\eta_{\widetilde{p}}}(z_{2})\right\}$\\
$\underline{\eta_{p}}(z_{1},z_{2})=\bigvee\limits_{z_{3}\in[z_{1}]_{\omega}[z_{2}]_{\omega}}\eta_{p}(z_{3})$\\
$\mbox{\hspace{1.3cm}}\leq\bigvee\limits_{z_{3}\in[z_{1}z_{2}]_{\omega}}\eta_{p}(z_{3})$\\$\mbox{\hspace{1.3cm}}=\bigvee\limits_{MN\in[z_{1}]_{\omega}[z_{2}]_{\omega}}\eta_{p}(MN)$\\$\mbox{\hspace{1.3cm}}\leq\bigvee\limits_{MN\in[z_{1}z_{2}]_{\omega}}\eta_{p}(MN)$\\
$\mbox{\hspace{1.3cm}}=\left(\bigvee\limits_{M\in[z_{1}]_{\omega}}\eta_{p}(M)\right)\bigwedge\left(\bigvee\limits_{N\in[z_{2}]_{\omega}}\eta_{p}(N)\right)$\\implies $\underline{\eta_{p}}(z_{1}z_{2})\leq max\left\{\underline{\eta_{p}}(M),\underline{\eta_{p}}(N)\right\}$.\\
\end{proof}
\end{theorem}
\begin{theorem}
Let $\omega$ be a $CR_{\omega}$ on $S$, and $P^{\Box}$ is a $P_{LI}^{\Box}$ (resp. $P_{RI}^{\Box}$) of $S$. Then $\underline{App}(P^{\Box})$ is a $P_{LI}^{\Box}$ (resp. $P_{RI}^{\Box}$) of $S$.
\begin{proof}
Since $\omega$ is a $CR_{\omega}$ on $S$, we have for all $z_{1},z_{2}\in S$ it follows that $[z_{1}][z_{2}]\subseteq[z_{1}z_{2}]_{\omega}$. We need to show that $\underline{App}(P^{\Box})=\left(\underline{\zeta_{p^{\Box}}},\underline{\eta_{p^{\Box}}}\right)=\left\langle\left[\underline{\zeta_{\widetilde{p}}},\underline{\eta_{\widetilde{p}}}\right],\left(\underline{\zeta_{p}},\underline{\eta_{p}}\right)\right\rangle$ is a $P_{LI}^{\Box}$ of $S$.\\Consider\\
$\underline{\zeta_{\widetilde{p}}}(z_{1}z_{2})=\bigwedge\limits_{z_{3}\in[z_{1}z_{2}]_{\omega}}\zeta_{\widetilde{p}}(z_{3})\geq\bigwedge\limits_{z_{3}\in[z_{1}]_{\omega}[z_{2}]_{\omega}}\zeta_{\widetilde{p}}(z_{3})$\\$\mbox{\hspace{1.3cm}}=\bigwedge\limits_{MN\in[z_{1}]_{\omega}[z_{2}]_{\omega}}\zeta_{\widetilde{p}}(MN)\geq\bigwedge\limits_{N\in[z_{2}]_{\omega}}\zeta_{\widetilde{p}}(N)$\\
implies $\underline{\zeta_{\widetilde{p}}}(z_{1}z_{2})\geq\underline{\zeta_{\widetilde{p}}}(z_{2})$\\
$\underline{\zeta_{p}}(z_{1}z_{2})=\bigwedge\limits_{z_{3}\in[z_{1}z_{2}]_{\omega}}\zeta_{p}(z_{3})\leq\bigwedge\limits_{z_{3}\in[z_{1}]_{\omega}[z_{2}]_{\omega}}\zeta_{p}(z_{3})$\\$\mbox{\hspace{1.3cm}}=\bigwedge\limits_{MN\in[z_{1}]_{\omega}[z_{2}]_{\omega}}\zeta_{p}(MN)\leq\bigwedge\limits_{N\in[z_{2}]_{\omega}}\zeta_{p}(N)$\\implies $\underline{\zeta_{p}}(z_{1}z_{2})\leq\underline{\zeta_{p}}(z_{2})$\\Next\\
$\underline{\eta_{\widetilde{p}}}(z_{1}z_{2})=\bigvee\limits_{z_{3}\in[z_{1}z_{2}]_{\omega}}\eta_{\widetilde{p}}(z_{3})\geq\bigvee\limits_{z_{3}\in[z_{1}]_{\omega}[z_{2}]_{\omega}}\eta_{\widetilde{p}}(z_{3})$\\$\mbox{\hspace{1.3cm}}=\bigvee\limits_{MN\in[z_{1}]_{\omega}[z_{2}]_{\omega}}\eta_{\widetilde{p}}(MN)\geq\bigvee\limits_{N\in[z_{2}]_{\omega}}\eta_{\widetilde{p}}(N)$\\
implies $\underline{\eta_{\widetilde{p}}}(z_{1}z_{2})\geq\underline{\eta_{\widetilde{p}}}(z_{2})$\\
$\underline{\eta_{p}}(z_{1}z_{2})=\bigvee\limits_{z_{3}\in[z_{1}z_{2}]_{\omega}}\eta_{p}(z_{3})\leq\bigvee\limits_{z_{3}\in[z_{1}]_{\omega}[z_{2}]_{\omega}}\eta_{p}(z_{3})$\\$\mbox{\hspace{1.3cm}}\bigvee\limits_{MN\in[z_{1}]_{\omega}[z_{2}]_{\omega}}\eta_{p}(MN)\leq\bigvee\limits_{N\in[z_{2}]_{\omega}}\eta_{p}(N)$\\implies $\underline{\eta_{p}}(z_{1}z_{2})\leq\underline{\eta_{p}}(z_{2})$\\ implies that $\underline{App}(P^{\Box})$ is a $\underline{P_{LI}^{\Box}}$ of $S$. Similarly, $\underline{App}(P^{\Box})$ is a $\underline{P_{RI}^{\Box}}$ of $S$.
\end{proof}
\end{theorem}
\begin{theorem}
Let $\omega$ be a $CR_{\omega}$ on semigroup $S$. If $P^{\Box}$ is a $P_{BI}^{\Box}$ of $S$. Then $\overline{App}(P^{\Box})$ is a $P_{BI}^{\Box}$ of $S$.
\begin{proof}
Since $\omega$ is a $CR_{\omega}$ on the semigroup $S$, we have for all $z_{1},z_{2},z_{3}\in S$\\ $[z_{1}]_{\omega}[z_{2}]_{\omega}[z_{3}]_{\omega}\subseteq[z_{1}z_{2}z_{3}]_{\omega}$, now show that $\overline{App}(P^{\Box})=\left(\overline{\zeta_{p^{\Box}}},\overline{\eta_{p^{\Box}}}\right)$ is a $P_{BI}^{\Box}$ of $S$. Consider the following\\
$\overline{\zeta_{\widetilde{p}}}(z_{1}z_{2}z_{3})=\bigvee\limits_{z\in[z_{1}z_{2}z_{3}]_{\omega}}\zeta_{\widetilde{p}}(z)\geq\bigvee\limits_{z\in[z_{1}]_{\omega}[z_{2}]_{\omega}[z_{3}]_{\omega}}\zeta_{\widetilde{p}}(z)$\\$\mbox{\hspace{1.3cm}}=\bigvee\limits_{abc\in[z_{1}]_{\omega}[z_{2}]_{\omega}[z_{3}]_{\omega}}\zeta_{\widetilde{p}}(abc)=\bigvee\limits_{a\in[z_{1}]_{\omega}b\in[z_{2}]_{\omega}c\in[z_{3}]_{\omega}}\zeta_{\widetilde{p}}(abc)$\\
$\mbox{\hspace{1.3cm}}\geq\bigvee\limits_{a\in[z_{1}]_{\omega}c\in[z_{3}]_{\omega}}\left\{\zeta_{\widetilde{p}}(a)\bigwedge\zeta_{\widetilde{p}}(c)\right\}$\\$\mbox{\hspace{1.3cm}}=\left\{\bigvee\limits_{a\in[z_{1}]_{\omega}}\zeta_{\widetilde{p}}(a)\right\}\bigwedge\left\{\bigvee\limits_{c\in[z_{3}]_{\omega}}\zeta_{\widetilde{p}}(c)\right\}$\\implies $\overline{\zeta_{\widetilde{p}}}(z_{1}z_{2}z_{3})\geq min\left\{\overline{\zeta_{\widetilde{p}}}(z_{1}),\overline{\zeta_{\widetilde{p}}}(z_{3})\right\}$\\
$\overline{\zeta_{p}}(z_{1}z_{2}z_{3})=\bigvee\limits_{z\in[z_{1}z_{2}z_{3}]_{\omega}}\zeta_{\widetilde{p}}(z)\leq\bigvee\limits_{z\in[z_{1}]_{\omega}[z_{2}]_{\omega}[z_{3}]_{\omega}}\zeta_{p}(z)$\\$\mbox{\hspace{1.3cm}}=\bigvee\limits_{abc\in[z_{1}]_{\omega}[z_{2}]_{\omega}[z_{3}]_{\omega}}\zeta_{p}(abc)=\bigvee\limits_{a\in[z_{1}]_{\omega}b\in[z_{2}]_{\omega}c\in[z_{3}]_{\omega}}\zeta_{p}(abc)$\\
$\mbox{\hspace{1.3cm}}\leq\bigvee\limits_{a\in[z_{1}]_{\omega}c\in[z_{3}]_{\omega}}\left\{\zeta_{p}(a)\bigvee\zeta_{p}(c)\right\}$\\
$\mbox{\hspace{1.3cm}}=\left\{\bigvee\limits_{a\in[z_{1}]_{\omega}}\zeta_{p}(a)\right\}\bigvee\left\{\bigvee\limits_{c\in[z_{3}]_{\omega}}\zeta_{p}(c)\right\}$\\implies $\overline{\zeta_{p}}(z_{1}z_{2}z_{3})\leq max\left\{\overline{\zeta_{p}}(z_{1}),\overline{\zeta_{p}}(z_{3})\right\}$\\
Next\\
$\overline{\eta_{\widetilde{p}}}(z_{1}z_{2}z_{3})=\bigwedge\limits_{z\in[z_{1}z_{2}z_{3}]_{\omega}}\eta_{\widetilde{p}}(z)\geq\bigwedge\limits_{z\in[z_{1}]_{\omega}[z_{2}]_{\omega}[z_{3}]_{\omega}}\eta_{\widetilde{p}}(z)$\\$\mbox{\hspace{1.3cm}}=\bigwedge\limits_{abc\in[z_{1}]_{\omega}[z_{2}]_{\omega}[z_{3}]_{\omega}}\eta_{\widetilde{p}}(abc)=\bigwedge\limits_{a\in[z_{1}]_{\omega}b\in[z_{2}]_{\omega}c\in[z_{3}]_{\omega}}\eta_{\widetilde{p}}(abc)$\\
$\mbox{\hspace{1.3cm}}\geq\bigwedge\limits_{a\in[z_{1}]_{\omega}c\in[z_{3}]_{\omega}}\left\{\eta_{\widetilde{p}}(a)\bigwedge\eta_{\widetilde{p}}(c)\right\}$\\$\mbox{\hspace{1.3cm}}=\left\{\bigwedge\limits_{a\in[z_{1}]_{\omega}}\eta_{\widetilde{p}}(a)\right\}\bigwedge\left\{\bigwedge\limits_{c\in[z_{3}]_{\omega}}\eta_{\widetilde{p}}(c)\right\}$\\implies $\overline{\eta_{\widetilde{p}}}(z_{1}z_{2}z_{3})\geq min\left\{\overline{\eta_{\widetilde{p}}}(z_{1}),\overline{\eta_{\widetilde{p}}}(z_{3})\right\}$\\
$\overline{\eta_{p}}(z_{1}z_{2}z_{3})=\bigwedge\limits_{z\in[z_{1}z_{2}z_{3}]_{\omega}}\eta_{\widetilde{p}}(z)\leq\bigwedge\limits_{z\in[z_{1}]_{\omega}[z_{2}]_{\omega}[z_{3}]_{\omega}}\eta_{p}(z)$\\$\mbox{\hspace{1.3cm}}=\bigwedge\limits_{abc\in[z_{1}]_{\omega}[z_{2}]_{\omega}[z_{3}]_{\omega}}\eta_{p}(abc)=\bigwedge\limits_{a\in[z_{1}]_{\omega}b\in[z_{2}]_{\omega}c\in[z_{3}]_{\omega}}\eta_{p}(abc)$\\
$\mbox{\hspace{1.3cm}}\leq\bigwedge\limits_{a\in[z_{1}]_{\omega}c\in[z_{3}]_{\omega}}\left\{\eta_{p}(a)\bigvee\eta_{p}(c)\right\}$\\
$\mbox{\hspace{1.3cm}}=\left\{\bigwedge\limits_{a\in[z_{1}]_{\omega}}\eta_{p}(a)\right\}\bigvee\left\{\bigwedge\limits_{c\in[z_{3}]_{\omega}}\eta_{p}(c)\right\}$\\implies $\overline{\eta_{p}}(z_{1}z_{2}z_{3})\leq max\left\{\overline{\eta_{p}}(z_{1}),\overline{\eta_{p}}(z_{3})\right\}$
\end{proof}
\end{theorem}
\begin{theorem}
Let $\omega$ be a complete $CR_{\omega}$ on semigroup $S$. Let $P^{\Box}$ is a $P_{BI}^{\Box}$ of $S$. Then $\underline{App}(P^{\Box})$ is a $P_{BI}^{\Box}$ of $S$.
\begin{proof}
Since $\omega$ is a $CR_{\omega}$ on the semigroup $S$, we have for all $z_{1},z_{2},z_{3}\in S$\\ $[z_{1}]_{\omega}[z_{2}]_{\omega}[z_{3}]_{\omega}\subseteq[z_{1}z_{2}z_{3}]_{\omega}$, we show that $\underline{App}(P^{\Box})=\left(\underline{\zeta_{p^{\Box}}},\underline{\eta_{p^{\Box}}}\right)$ is a $P_{BI}^{\Box}$ of $S$.\\ Consider the following\\
$\underline{\zeta_{\widetilde{p}}}(z_{1}z_{2}z_{3})=\bigwedge\limits_{z\in[z_{1}z_{2}z_{3}]_{\omega}}\zeta_{\widetilde{p}}(z)\geq\bigwedge\limits_{z\in[z_{1}]_{\omega}[z_{2}]_{\omega}[z_{3}]_{\omega}}\zeta_{\widetilde{p}}(z)$\\$\mbox{\hspace{1.3cm}}=\bigwedge\limits_{abc\in[z_{1}]_{\omega}[z_{2}]_{\omega}[z_{3}]_{\omega}}\zeta_{\widetilde{p}}(abc)=\bigwedge\limits_{a\in[z_{1}]_{\omega}b\in[z_{2}]_{\omega}c\in[z_{3}]_{\omega}}\zeta_{\widetilde{p}}(abc)$\\
$\mbox{\hspace{1.3cm}}\geq\bigwedge\limits_{a\in[z_{1}]_{\omega}c\in[z_{3}]_{\omega}}\left\{\zeta_{\widetilde{p}}(a)\bigwedge\zeta_{\widetilde{p}}(c)\right\}$\\$\mbox{\hspace{1.3cm}}=\left\{\bigwedge\limits_{a\in[z_{1}]_{\omega}}\zeta_{\widetilde{p}}(a)\right\}\bigwedge\left\{\bigwedge\limits_{c\in[z_{3}]_{\omega}}\zeta_{\widetilde{p}}(c)\right\}$\\implies $\underline{\zeta_{\widetilde{p}}}(z_{1}z_{2}z_{3})\geq min\left\{\underline{\zeta_{\widetilde{p}}}(z_{1}),\underline{\zeta_{\widetilde{p}}}(z_{3})\right\}$\\
$\underline{\zeta_{p}}(z_{1}z_{2}z_{3})=\bigwedge\limits_{z\in[z_{1}z_{2}z_{3}]_{\omega}}\zeta_{\widetilde{p}}(z)\leq\bigwedge\limits_{z\in[z_{1}]_{\omega}[z_{2}]_{\omega}[z_{3}]_{\omega}}\zeta_{p}(z)$\\$\mbox{\hspace{1.3cm}}=\bigwedge\limits_{abc\in[z_{1}]_{\omega}[z_{2}]_{\omega}[z_{3}]_{\omega}}\zeta_{p}(abc)=\bigvee\limits_{a\in[z_{1}]_{\omega}b\in[z_{2}]_{\omega}c\in[z_{3}]_{\omega}}\zeta_{p}(abc)$\\
$\mbox{\hspace{1.3cm}}\leq\bigwedge\limits_{a\in[z_{1}]_{\omega}c\in[z_{3}]_{\omega}}\left\{\zeta_{p}(a)\bigvee\zeta_{p}(c)\right\}$\\
$\mbox{\hspace{1.3cm}}=\left\{\bigwedge\limits_{a\in[z_{1}]_{\omega}}\zeta_{p}(a)\right\}\bigvee\left\{\bigwedge\limits_{c\in[z_{3}]_{\omega}}\zeta_{p}(c)\right\}$\\implies $\underline{\zeta_{p}}(z_{1}z_{2}z_{3})\leq max\left\{\underline{\zeta_{p}}(z_{1}),\underline{\zeta_{p}}(z_{3})\right\}$\\
Next\\
$\underline{\eta_{\widetilde{p}}}(z_{1}z_{2}z_{3})=\bigvee\limits_{z\in[z_{1}z_{2}z_{3}]_{\omega}}\eta_{\widetilde{p}}(z)\geq\bigvee\limits_{z\in[z_{1}]_{\omega}[z_{2}]_{\omega}[z_{3}]_{\omega}}\eta_{\widetilde{p}}(z)$\\$\mbox{\hspace{1.3cm}}=\bigvee\limits_{abc\in[z_{1}]_{\omega}[z_{2}]_{\omega}[z_{3}]_{\omega}}\eta_{\widetilde{p}}(abc)=\bigvee\limits_{a\in[z_{1}]_{\omega}b\in[z_{2}]_{\omega}c\in[z_{3}]_{\omega}}\eta_{\widetilde{p}}(abc)$\\
$\mbox{\hspace{1.3cm}}\geq\bigvee\limits_{a\in[z_{1}]_{\omega}c\in[z_{3}]_{\omega}}\left\{\eta_{\widetilde{p}}(a)\bigwedge\eta_{\widetilde{p}}(c)\right\}$\\$\mbox{\hspace{1.3cm}}=\left\{\bigvee\limits_{a\in[z_{1}]_{\omega}}\eta_{\widetilde{p}}(a)\right\}\bigwedge\left\{\bigvee\limits_{c\in[z_{3}]_{\omega}}\eta_{\widetilde{p}}(c)\right\}$\\implies $\underline{\eta_{\widetilde{p}}}(z_{1}z_{2}z_{3})\geq min\left\{\underline{\eta_{\widetilde{p}}}(z_{1}),\underline{\eta_{\widetilde{p}}}(z_{3})\right\}$\\
$\underline{\eta_{p}}(z_{1}z_{2}z_{3})=\bigvee\limits_{z\in[z_{1}z_{2}z_{3}]_{\omega}}\eta_{\widetilde{p}}(z)\leq\bigvee\limits_{z\in[z_{1}]_{\omega}[z_{2}]_{\omega}[z_{3}]_{\omega}}\eta_{p}(z)$\\$\mbox{\hspace{1.3cm}}=\bigvee\limits_{abc\in[z_{1}]_{\omega}[z_{2}]_{\omega}[z_{3}]_{\omega}}\eta_{p}(abc)=\bigvee\limits_{a\in[z_{1}]_{\omega}b\in[z_{2}]_{\omega}c\in[z_{3}]_{\omega}}\eta_{p}(abc)$\\
$\mbox{\hspace{1.3cm}}\leq\bigvee\limits_{a\in[z_{1}]_{\omega}c\in[z_{3}]_{\omega}}\left\{\eta_{p}(a)\bigvee\eta_{p}(c)\right\}$\\
$\mbox{\hspace{1.3cm}}=\left\{\bigvee\limits_{a\in[z_{1}]_{\omega}}\eta_{p}(a)\right\}\bigvee\left\{\bigvee\limits_{c\in[z_{3}]_{\omega}}\eta_{p}(c)\right\}$\\implies $\underline{\eta_{p}}(z_{1}z_{2}z_{3})\leq max\left\{\underline{\eta_{p}}(z_{1}),\underline{\eta_{p}}(z_{3})\right\}$
\end{proof}
\end{theorem}
\begin{theorem}
Let $\omega$ be a $CR_{\omega}$ on semigroup $S$. If $P^{\Box}$ is a $P_{II}^{\Box}$ of $S$. Then $\overline{App}(P^{\Box})$ is a $P_{II}^{\Box}$ of $S$.
\begin{proof}
Proof directly follow from theorem 4.5
\end{proof}
\end{theorem}
\begin{theorem}
Let $\omega$ be a complete $CR_{\omega}$ on semigroup $S$. Let $P^{\Box}$ is a $P_{II}^{\Box}$ of $S$. Then $\underline{App}(P^{\Box})$ is a $P_{II}^{\Box}$ of $S$.
\begin{proof}
Proof directly follow from theorem 4.6
\end{proof}
\end{theorem}
\section{Conclusions}
Cubic Pythagorean fuzzy sets are the generalization of cubic sets. In this paper, we have presented the concept of rough cubic Pythagorean fuzzy sets in semigroups, which can handle the vagueness in a proactive way than cubic sets. Then, we have extended the notion of rough cubic Pythagorean fuzzy sets to the lower and upper approximations of Pythagorean fuzzy left (right)ideals, bi-ideals, interior ideals in semigroups and also discussed some of its related properties. We aim to extend this work to some algebraic structures namely gamma semigroup, Po-gamma-semigroup, and subtraction semigroup.

%

\bigskip

\end{document}